\documentclass[11pt]{article}
\usepackage[margin=1in]{geometry}
\usepackage{amsmath,amssymb,amsthm,mathtools,latexsym,amssymb,parskip,enumerate}
\def\reals{\mathbb{R}}
\usepackage{hyperref}
\newtheorem{theorem}{Theorem}
\newtheorem{lemma}{Lemma}

\newcommand{\E}{\mathbb E}
\newcommand{\Pp}{\mathbb P}

\newtheorem{claim}{Claim}
\newtheorem{proposition}{Proposition}
\def\bp{{\bf p}}\def\bu{{\bf u}}\def\by{{\bf y}}\def\bz{{\bf z}}\def\bb{{\bf b}}\def\be{{\bf e}}\def\bx{{\bf x}}\def\bv{{\bf v}}
\def\bp{{\bf p}}\def\bc{{\bf c}}
\DeclareMathOperator{\Ent}{Ent}

\title{Contiguity of unions of perfect matchings\\and random regular uniform hypergraphs}
\author{Alan Frieze\thanks{Department of Mathematical Sciences, Carnegie Mellon University, Pittsburgh PA, USA, 15213. Research supported in part by NSF grant DMS1952285.}  and Xavier P\'erez-Gim\'enez\thanks{Department of Mathematics, University of Nebraska-Lincoln, Lincoln NE, USA, 68588. Research supported in part by Simons Foundation Grant \#587019 and by NSF grant DMS2201590.}} 
\date{}
\begin{document}
\maketitle

\begin{abstract}
For fixed integers $s\ge3$ and $r$ satisfying
$\log(r!)>\frac{r(s-1)}s\log r$, the union of $r$ independent uniform perfect matchings of the complete $s$-uniform hypergraph is mutually contiguous to the uniform simple $r$-regular $s$-uniform hypergraph, as $n\to\infty$ through multiples of $s$. In particular, $r=\lceil e^s\rceil$ suffices. The main input is the sharp entropy inequality for random permutations of Bristiel and Caputo. It proves the global maximization needed for the second moment of the number of complete factorizations. The remaining calculation, including its lattice factor, matches the short-cycle contribution required by small-subgraph conditioning.\\
This is an expanded version of a proof provided by ChatGPT-Astra.
\end{abstract}
\section{Introduction}
A random $r$-regular graph, $r\geq3$, on $n$ vertices is contiguous to the union of $r$ random perfect matchings provided $n$ is even, conditioned on the union of the matchings being simple. This has been known for some time, see for example Wormald \cite{Wreg}. In this paper we extend this decomposition to random regular hypergraphs.
\begin{theorem}\label{thm:main}
Fix integers $s\ge3$ and $r\ge2$ such that
\begin{equation}\label{eq:condition}
 f_{r,s}:=\log(r!)-\frac{r(s-1)}s\log r>0.
\end{equation}
Let $n\to\infty$ through multiples of $s$. Let $U_{n,r,s}$ be the multiset union of $r$ independent uniform perfect matchings of the complete $s$-uniform hypergraph on $[n]$. Let $H_{n,r,s}$ be uniform over all simple $r$-regular $s$-uniform hypergraphs on $[n]$. These two sequences are mutually contiguous: for every sequence of events $\mathcal E_n$,
\[
 \Pp(U_{n,r,s}\in\mathcal E_n)\longrightarrow0 \quad\Longleftrightarrow\quad
 \Pp(H_{n,r,s}\in\mathcal E_n)\longrightarrow0.
\]
The same conclusion holds with ordinary set union, since repeated edges occur with probability $o(1)$. 
\end{theorem}
Condition \eqref{eq:condition} implies $r>s$. Indeed, $r!\le r^{r-1}$, so
\[
 \frac{s-1}{s}<\frac{\log(r!)}{r\log r} \le\frac{r-1}{r}.
\]
Also,
\[
 \log(r!)>\int_1^r\log x\,dx=r\log r-r+1,
\]
so \eqref{eq:condition} holds whenever $\log r\ge s$.

One motivation for proving Theorem \ref{thm:main} comes from our work \cite{Loose} on the existence of loose Hamilton cycles in random hypergraphs. A loose Hamilton cycle in an $s$-uniform hypergraph is a sequence $\be_1,\be_2,\ldots,\be_m,\be_{m+1}=\be_1,m=n/(s-1)$ of edges such that (i) $|\be_i\cap \be_{i+1}|=1$ for $i=1,2,\ldots,m$, (ii) $\bigcup_{i=1}^m\be_i=[n]$ and (iii) non-consecutive edges are disjoint. In this paper we prove 
\begin{theorem}\label{th1}
The threshold for the binomial random hypergraph $H_{n,p;s}$, $(s-1)\mid n$,  to have a loose Hamilton cycle is at $p\sim \log n/\binom{n-1}{s-1}$. 
\end{theorem}
As part of the proof of this result, we showed that if $s\mid n$ and $p=(1+\epsilon)\log n/\binom{n-1}{s-1}$ then
\begin{equation}\label{weshowed}
 \text{w.h.p. $H_{n,p;s}$ contains $r$ independent uniform perfect matchings.}
\end{equation}
 Here $r$ can be made as large as needed. Now  Altman, Greenhill,  Isaev and Ramadurai \cite{AGIR} proved the following theorem:
\begin{theorem}\label{thm:lcr}
Suppose that $s\mid rn$ and $(s-1)\mid n$. Then,
\[
\lim_{n\to\infty}\mathbf P(H_{n,r,s}\text{ contains a loose Hamilton cycle})=\begin{cases}1&r>\rho(s).\\0&r\leq \rho(s).\end{cases}
\]
Here $\rho=\rho(s)$ is the unique real in $(2,\infty)$ such that
\[
(\rho-1)(s-1)\left(\frac{\rho s-\rho-s}{\rho s-\rho}\right)^{(s-1)(\rho s-\rho-s)/s}=1.
\]
\end{theorem}
Thus Theorem \ref{th1} follows from \eqref{weshowed}, Theorem  \ref{thm:main} and Theorem \ref{thm:lcr}, after dealing with the extra condition $s\mid n$ using an idea of Ferber \cite{Ferber} to show that the probability of loose Hamilton cycle is bounded below by a constant plus a result of Friedgut \cite{Fried}, showing that there are no course thresholds when the number of random edges $m=cn\log n$. (Note that when $p=(1-\epsilon)\log n/\binom{n-1}{s-1}$ then w.h.p. there will be isolated vertices in $H_{n,p;s}$, clearly ruling out a loose hamilton cycle.)
\paragraph{Use of AI} We used ChatGPT to come up with an initial proof. This expanded write-up makes the proof more or less self-contained.
\section{Configuration model and the change of measure}

Give each vertex $v$, $r$ distinguishable points $W_v$ and partition all $rn$ points $W=\bigcup_{v\in[n]}W_v$ uniformly into $rn/s$ sets that represent hyperedges of size $s$.  (We will loosely call these sets hyperedges.) A hyperedge containing two points at the same vertex is called a loop. Write $P_n$ for this configuration law. For nonnegative integers $m$ divisible by $s$, put 
\[
 T_s(m)=\frac{m!}{(m/s)!\,(s!)^{m/s}},\qquad T_s(0)=1,
\]
and set $T_s(m)=0$ if $s\nmid m$.

Let $Z$ count {\em valid} colourings of the hyperedges with colours $[r]$ such that the $r$ points at every vertex have distinct colours. Here the color of a point is the color of the hyperedge containing it. Thus each colour class projects to a perfect matching. In particular, $Z=0$ whenever the configuration has a loop. Thus a valid colouring is exactly an \emph{ordered factorisation}
\[
F=F_1\uplus F_2\uplus\cdots\uplus F_r,
\]
where each $F_i$ is a perfect matching and the union is over hyperedge occurrences. Therefore, $Z(F)$ counts these labelled factorisations.
If factorisations are regarded as unordered, each one is counted $r!$ times by $Z(F)$.

At each vertex there are $r!$ bijections from its points to $[r]$. For a fixed choice of these bijections, there are $T_s(n)^r$ partitions into monochromatic hyperedges. ($T_s(n)$ counts the number of choices for making perfect matchings from the points.) Hence \begin{equation}\label{eq:firstmoment}
 \E_{P_n}Z=(r!)^n\frac{T_s(n)^r}{T_s(rn)}.
\end{equation}
We next define a probability distribution $Q_n$ on configurations. For every configuration $F$, we let 
\begin{equation}\label{eq:Qexplicit}
 Q_n(F) =\frac{Z(F)P_n(F)}{\E_{P_n}Z} =\frac{Z(F)}{(r!)^nT_s(n)^r}.
\end{equation}
Equivalently, sample $Q_n$ as follows. At each vertex, independently choose a uniform bijection from its $r$ points to the colours $[r]$.
For each colour, independently partition its $n$ points uniformly into hyperedges of size $s$. Finally, forget the colours. Each valid coloured configuration has probability $((r!)^nT_s(n)^r)^{-1}$, and an uncoloured configuration $F$ has exactly $Z(F)$ such colourings, proving \eqref{eq:Qexplicit}. After replacing each point by the vertex to which it belongs this construction gives the union of $r$ independent uniform perfect matchings. 

We first prove mutual contiguity between $Q_n$ and the configuration model conditioned to have no loops. At the end we show that simplicity has probability $1-o(1)$ under both laws. The configuration model conditioned on simplicity, after projection to $[n]$, gives exactly the uniform simple $r$-regular $s$-uniform hypergraph. These facts yield the main theorem.
\section{The second moment}
Let $u$ denote the uniform probability distribution on $S_r$, the set of permutations of $[r]$. For a probability distribution $\bp=(\bp_\sigma)_{\sigma\in S_r}$, write
\[
 a_{ij}=\sum_{\sigma:\,\sigma(i)=j}p_\sigma, \qquad D(\bp\Vert \bu)=\sum_{\sigma}p_\sigma\log(r!p_\sigma),
\]
and put
\[
 B(\bp)=\sum_{i,j}a_{ij}\log(ra_{ij}).
\]
Inequality (4.13) of Bristiel and Caputo \cite{BC} can be used to prove
\begin{equation}\label{eq:BC}
 B(\bp)\le \frac{r\log r}{\log(r!)}D(\bp\Vert \bu).
\end{equation}
Here \(0\log0=0\). For a proof of \eqref{eq:BC}, using \cite{BC}, see Sections \ref{sec:BC}, \ref{BC1} in the appendix.

Now let $\sigma=(\sigma_1,\ldots,\sigma_n)$ be independent uniform permutations of $[r]$, and let $P_{\sigma_v},v\in [n]$ be their permutation matrices. For the sampled permutations, let
\begin{align}
p_\rho&=\frac1n\bigl|\{v\in[n]:\sigma_v=\rho\}\bigr|, \qquad  a_{ij}:=\sum_{\rho:\rho(i)=j}p_\rho.\nonumber\\
 N&=\sum_{v=1}^n P_{\sigma_v}, \qquad N_{ij}=|\{v:\sigma_v(i)=j\}|=na_{ij}.\label{N=}
\end{align}
Every row and column of $N$ sums to $n$.

Fix the point-to-colour assignments for a first factorization. At vertex $v$, the assignment for a second factorization is specified by a permutation $\sigma_v$ of the first colours. A \emph{cell} is the set of points that have a specified pair of colours: colour $i$ in the first factorization and colour $j$ in the second. At vertex $v$, the point with first colour $i$ has second colour $\sigma_v(i)$, so it belongs to cell $(i,\sigma_v(i))$.

We first prove
\begin{lemma}\label{cl1}
\begin{equation}\label{eq:secondexact}
 \frac{\E_{P_n}Z^2}{(\E_{P_n}Z)^2} = \frac{T_s(rn)}{T_s(n)^{2r}}\, \E_{\sigma_1,\ldots,\sigma_n} \left[\prod_{i,j}T_s(N_{ij})\right].
\end{equation}
\end{lemma}
\begin{proof}
For each configuration $F$, the quantity $Z(F)^2$ counts ordered pairs of valid colourings. Since there are $T_s(rn)$ equally likely configurations,
\[
  \mathbb E_{P_n}Z^2  =\frac{1}{T_s(rn)}\sum_F Z(F)^2.
\]
Count these pairs by first fixing the first colouring. There are $(r!)^n$ choices of its point-to-colour assignments, one bijection at each vertex. To specify the second colouring, choose at each vertex $v$ a permutation $\sigma_v\in S_r$: the point with first colour $i$ gets second colour $\sigma_v(i)$. Thus there are $(r!)^n$ possible lists $(\sigma_1,\ldots,\sigma_n)$.

For a fixed list of permutations, there are  $N_{ij}$ points in cell $(i,j)$ whose first colour is $i$ and whose second colour is $j$. A configuration compatible with both colourings must partition each cell separately into $s$-element hyperedges. There are $T_s(N_{ij})$ ways to partition cell $(i,j)$, so the number of compatible configurations is $\prod_{i,j}T_s(N_{ij})$. Therefore
\[
\mathbb E_{P_n}Z^2  =\frac{(r!)^n}{T_s(rn)}\sum_{\sigma_1,\ldots,\sigma_n\in S_r}\prod_{i,j}T_s(N_{ij}).
\]
Also, 
\[
  \mathbb E_{P_n}Z=\frac{(r!)^nT_s(n)^r}{T_s(rn)}.
\]
Dividing the expression for $\mathbb E_{P_n}Z^2$ by $(\mathbb E_{P_n}Z)^2$ gives
\[
\frac{\mathbb E_{P_n}Z^2}{(\mathbb E_{P_n}Z)^2}=\frac{T_s(rn)}{T_s(n)^{2r}}\,
\frac{1}{(r!)^n}\sum_{\sigma_1,\ldots,\sigma_n\in S_r}\prod_{i,j}T_s(N_{ij}).
\]
The final factor is the expectation when $\sigma_1,\ldots,\sigma_n$ are chosen independently and uniformly from $S_r$. This is exactly equation \eqref{eq:secondexact}.
\end{proof}
Given \eqref{eq:secondexact}, we can prove
\begin{lemma}\label{lem:second}
Under \eqref{eq:condition}, with $d=(r-1)^2$,
\begin{equation}\label{eq:secondlimit}
 \frac{\E_{P_n}Z^2}{(\E_{P_n}Z)^2}\longrightarrow\left(\frac{r-1}{r-s}\right)^{d/2}.
\end{equation}
\end{lemma}
\begin{proof}
Fix $r>s$ satisfying
\[
  \log(r!)>\frac{r(s-1)}s\log r,
\]
and let $n\to\infty$ through multiples of $s$. Let
\[
K=r!,\qquad q=K-1,\qquad d=(r-1)^2,\qquad  \alpha=\frac{s-1}{s}.
\]
For a vector of counts $\bc=(c_\sigma)_{\sigma\in S_r}$ with $\sum_\sigma c_\sigma=n$, let
\[
  p_\sigma=\frac{c_\sigma}{n},\qquad N=\sum_\sigma c_\sigma P_\sigma, 
\]
where $(P_\sigma)_{ij}=\mathbf 1_{\{\sigma(i)=j\}}$. Define
\[
  \mu_n(\bc)=\frac{n!}{K^n\prod_\sigma c_\sigma!},
  \qquad
  W_n(\bc)=\frac{T_s(rn)}{T_s(n)^{2r}}\prod_{i,j}T_s(N_{ij}).
\]
Equation \eqref{eq:secondexact} becomes 
\begin{equation}\label{Rn}
  R_n:=\frac{\mathbb E_{P_n}Z^2}{(\mathbb E_{P_n}Z)^2}= \sum_{\substack{c_\sigma\ge0\\\sum_\sigma c_\sigma=n}}\mu_n(\bc)W_n(\bc).
\end{equation}
Only counts satisfying $N_{ij}\equiv0\pmod s$ contribute.

The inequality \eqref{eq:BC} gives
\[
  -D(\bp\Vert \bu)+\alpha B(\bp)\le-\eta D(\bp\Vert \bu),  \qquad  \eta=1-\alpha\frac{r\log r}{\log(r!)}>0.
\]
\begin{claim}\label{Phi}
There exists $C=C_{r,s}$ such that for a fixed $\bc$ 
 \begin{equation}\label{BD}
 \mu_n(\bc)W_n(\bc)\leq n^{C}\exp\{n[-D(\bp\Vert \bu)+\alpha B(\bp)]\},
\end{equation}
where $\bp=\bc/n$.
\end{claim}
{\em Proof of Claim \ref{Phi}}
 For a fixed $\bc$, 
\begin{equation}\label{boundtype}
  \mu_n(\bc)=\frac{n!}{K^n\prod_{\sigma\in S_r}c_\sigma!} \le n^{C_r}\exp\{-nD(\bp\Vert \bu)\},
\end{equation}
where the RHS is obtained by using Stirling's formula.

Stirling's formula gives, for multiples $m$ of $s$,
\[
  \log T_s(m)=\alpha m(\log m-1)-\frac{m}{s}\log((s-1)!)+O_s(\log(m+1)).
 \]
Substitute this into $\log W_n(\bc)$. Since $\sum_{i,j}a_{ij}=r$, the terms linear in $m$ cancel, and the remaining terms of order $n$ are
\[
\begin{aligned}
\log W_n(\bc)&=\alpha n\left(r\log r+\sum_{i,j}a_{ij}\log a_{ij}\right)+O_{r,s}(\log n)\\
&=\alpha n\sum_{i,j}a_{ij}\log(ra_{ij})+O_{r,s}(\log n)\\
&=\alpha n B(\bp)+O_{r,s}(\log n).
\end{aligned}
\]
Therefore
\begin{equation}\label{eq:global}
  \mu_n(\bc)\,W_n(\bc)\le n^{C_{r,s}}\exp\!\left\{n\bigl[-D(\bp\Vert \bu)+\alpha B(\bp)\bigr]\right\}\leq n^{C_{r,s}}\exp\!\left\{-\eta nD(\bp\Vert \bu)\right\},
\end{equation}
where $\eta=1- \frac{\alpha r\log r}{\log(r!)}>0$.

The factor $n^{C_{r,s}}$ absorbs the polynomial factors from the two Stirling estimates. \\
{\em End of proof of Claim \ref{Phi}}

Pinsker's inequality implies that $D(\bp\Vert \bu)\ge\frac12\|p-u\|_1^2$. Since there are at most $(n+1)^K$ types, the total contribution
from $\|\bp-\bu\|_1>n^{-2/5}$ is
\begin{equation}\label{outside}
  O\!\left(n^{C+K}e^{-\eta n^{1/5}/2}\right)=o(1).
\end{equation}

For the remaining choices of $\bc$, put
\[
  y_\sigma=\frac{c_\sigma-n/K}{\sqrt n}.
\]
Then $y$ belongs to the $q$-dimensional Euclidean space
\begin{equation}\label{defH}
H=\left\{y\in\mathbb R^{S_r}:\sum_\sigma y_\sigma=0\right\},
\end{equation}
and $\|y\|_1\le n^{1/10}$. 

Define the linear map $L$ by 
\[
  L\by=\sum_{\sigma\in S_r}P_\sigma y_\sigma.
\]
This maps $H$ into
\begin{equation}\label{defV}
  V=\{M\in\mathbb R^{r\times r}: M\mathbf1=0,\ M^{\mathsf T}\mathbf1=0\},  \qquad \dim V=d,
\end{equation}
and note that
\[
 N=\frac nrJ+\sqrt n\,L\by,
\]
where \(J=\mathbf1\mathbf1^{\mathsf T}\), the \(r\times r\) all-ones matrix.
\begin{claim}\label{lem1}
If  \(\mathbf c\) is admissible, meaning \(s\mid N_{ij}\) for all \(i,j\), and $\|\bp-\bu\|_1\le n^{-2/5}$ then
\[
 \mu_n(\bc)W_n(\bc) = \frac{s^{d/2}K^{K/2}}{(2\pi n)^{q/2}}\exp\!\left\{-\frac K2\|\by\|^2+\frac{\alpha r}{2}\|L\by\|_F^2  \right\}  \bigl(1+O(n^{-1/5})\bigr).
\]
\end{claim}
{\em Proof of Claim \ref{lem1}}\\
Since $p_\sigma=1/K+O(n^{-2/5})$, every $c_\sigma\ge n/(2K)$ for sufficiently large $n$. Stirling's formula gives
\[
\begin{aligned}
  \mu_n(\bc)&=\frac{n!}{K^n\prod_\sigma c_\sigma!}=\frac{\sqrt{2\pi n}}{\prod_\sigma\sqrt{2\pi np_\sigma}}\,  \frac{(n/e)^n}       {K^n\prod_\sigma(np_\sigma/e)^{np_\sigma}} \bigl(1+O(n^{-1})\bigr)\\
  &=  (2\pi n)^{-(K-1)/2}\prod_\sigma p_\sigma^{-1/2}\exp\!\left\{-n\log K-n\sum_\sigma p_\sigma\log p_\sigma \right\}  \bigl(1+O(n^{-1})\bigr)\\
  &=(2\pi n)^{-q/2}\prod_\sigma p_\sigma^{-1/2} e^{-nD(\bp\Vert \bu)}
  \bigl(1+O(n^{-1})\bigr).
\end{aligned}
\]

To evaluate the product, put $\delta_\sigma=p_\sigma-1/K$. Then $\sum_\sigma\delta_\sigma=0$, and
\[
\begin{aligned}
\prod_\sigma p_\sigma^{-1/2}&=\exp\left\{\frac K2\log K  -\frac12\sum_\sigma\log(1+K\delta_\sigma)\right\}\\
&=\exp\left\{\frac K2\log K-\frac K2\sum_\sigma\delta_\sigma  +O\!\left(\sum_\sigma\delta_\sigma^2\right)\right\}\\
&=K^{K/2}\bigl(1+O(n^{-4/5})\bigr)
\end{aligned}
\]
since
\[
  \sum_\sigma\delta_\sigma^2\le\left(\sum_\sigma|\delta_\sigma|\right)^2  \le n^{-4/5}.
\]
\begin{proposition}\label{cl2}
\[
  nD(\bp\Vert \bu)=\frac K2\|\by\|^2+O_K(n^{-1/5}).
\]
\end{proposition}
{\em Proof of Proposition \ref{cl2}}

This follows from Taylor expansion about the uniform distribution $u_\sigma=1/K$. For $\phi(t)=t\log(Kt)$, we have
\[
  \phi(1/K)=0,\qquad  \phi'(1/K)=1,\qquad  \phi''(1/K)=K,\qquad  \phi'''(t)=-\frac1{t^2}.
\]
For sufficiently large $n$, every $p_\sigma$ lies between $1/(2K)$ and $3/(2K)$, so the third derivative is uniformly bounded there. Hence
\[
  \phi(1/K+\delta_\sigma)  =\delta_\sigma+\frac K2\delta_\sigma^2   +O_K(|\delta_\sigma|^3).
\]
Summing and using $\sum_\sigma\delta_\sigma=0$ gives
\[
  D(\bp\Vert \bu)  =\frac K2\sum_\sigma\delta_\sigma^2+O_K\!\left(\sum_\sigma|\delta_\sigma|^3\right).
\]
After multiplying by $n$, the quadratic term becomes
\[
  \frac{nK}{2}\sum_\sigma\delta_\sigma^2=\frac K2\|\by\|^2.
\]
The remainder satisfies
\[
  n\sum_\sigma|\delta_\sigma|^3\le n\left(\sum_\sigma|\delta_\sigma|\right)^3\le n(n^{-2/5})^3=n^{-1/5}.
\]
\noindent{\em End of proof of Proposition \ref{cl2}}

This gives
\[
  \mu_n(\bc)  =  \frac{K^{K/2}}{(2\pi n)^{q/2}}e^{-K\|\by\|^2/2}\bigl(1+O(n^{-1/5})\bigr).
\]
Next consider
\[
  W_n(\bc)  =\frac{T_s(rn)}{T_s(n)^{2r}}\prod_{i,j}T_s(N_{ij}).
\]
For admissible counts, all $N_{ij}$ are divisible by $s$. Put
\[
  F(t)=t(\log t-1),\qquad  \Delta_{ij}=N_{ij}-\frac nr.
\]
Since $\sum_\sigma P_\sigma=(K/r)J$, we have
\[
  N=\frac nrJ+\sqrt n\,L\by,  \qquad  \Delta_{ij}=\sqrt n\,(L\by)_{ij}.
\]
Moreover,
\[
  |\Delta_{ij}|\le n^{3/5},  \qquad  \sum_{i,j}\Delta_{ij}=0.
\]

Stirling's formula for $T_s$ is
\[
  \log T_s(m)  =\frac12\log s+\alpha F(m)   -\frac ms\log((s-1)!)+O(m^{-1}).
\]
Substitution into $W_n(\bc)$ gives
\[
  \log W_n(\bc)  =  \frac d2\log s  +\alpha\left[      F(rn)-2rF(n)+\sum_{i,j}F(N_{ij})    \right]+O(n^{-1}).
\]
Taylor expansion at $x=n/r$ gives
\[
  F(x+\Delta)
  =
  F(x)+\Delta\log x+\frac{r}{2n}\Delta^2+O(n^{-2}|\Delta|^3).
\]
The linear terms vanish upon summation, and
\[
  F(rn)-2rF(n)+r^2F(n/r)=0.
\]
Consequently,
\[
\begin{aligned}
  \log W_n(\bc)  &= \frac d2\log s+\frac{\alpha r}{2n}\sum_{i,j}\Delta_{ij}^2+O\!\left(n^{-2}\sum_{i,j}|\Delta_{ij}|^3+n^{-1}\right)\\
  &=\frac d2\log s+\frac{\alpha r}{2}\|L\by\|_F^2+O(n^{-1/5}),
\end{aligned}
\]
where the last error bound follows from $n^{-2}|\Delta_{ij}|^3\le n^{-1/5}$. Exponentiating,
\[
  W_n(\bc)=s^{d/2}e^{\alpha r\|L\by\|_F^2/2}\bigl(1+O(n^{-1/5})\bigr).
\]
Finally, multiplication gives
\[
  \mu_n(\bc)W_n(\bc) = \frac{s^{d/2}K^{K/2}}{(2\pi n)^{q/2}}\exp\!\left\{-\frac K2\|\by\|^2+\frac{\alpha r}{2}\|L\by\|_F^2  \right\}  \bigl(1+O(n^{-1/5})\bigr).
\]
{\em End of proof of Claim \ref{lem1}}

Taking $M=L\by$ in \eqref{defV} gives $\|L\by\|_F^2=\langle \by,L^TL\by\rangle_H$. Thus, with $I_H$ denoting the identity operator on $H$,
\[
  -\frac K2\|\by\|^2+\frac{\alpha r}{2}\|L\by\|_F^2 =-\frac12\langle \by,(KI_H-\alpha rL^TL)\by\rangle_H.
\]
Writing $Q=KI_H-\alpha rL^TL$ therefore yields 
\begin{equation}\label{muwq}
  \mu_n(\bc)W_n(\bc)  =\frac{s^{d/2}K^{K/2}}{(2\pi n)^{q/2}}e^{-\langle \by,Q\by\rangle/2}\bigl(1+O(n^{-1/5})\bigr).
\end{equation}
\begin{claim}\label{cl3}
\[
 LL^{\mathsf T}=\frac{K}{r-1}I_V.
\]
\end{claim}
{\em Proof of Claim \ref{cl3}}

If $B\in V$ then
\[
  (L^{\mathsf T}B)_\sigma=\sum_{k=1}^r B_{k,\sigma(k)}.
\]
For fixed $i,j$,
\[
  (LL^{\mathsf T}B)_{ij}=\sum_{\substack{\sigma\in S_r\\\sigma(i)=j}}\sum_{k=1}^r B_{k,\sigma(k)}=(r-1)!\,B_{ij}  +(r-2)!\sum_{\substack{k\ne i\\\ell\ne j}}B_{k\ell}.
\]
Indeed, there are $(r-1)!$ permutations with $\sigma(i)=j$. For $k\ne i$ and $\ell\ne j$, there are $(r-2)!$ permutations satisfying both $\sigma(i)=j$ and $\sigma(k)=\ell$. Now
\[
  \sum_{\substack{k\ne i\\\ell\ne j}}B_{k\ell}  =  \sum_{k,\ell}B_{k\ell}  -\sum_\ell B_{i\ell}  -\sum_k B_{kj}  +B_{ij}  =B_{ij}.
\]
Therefore
\[
  (LL^{\mathsf T}B)_{ij}
  =\bigl((r-1)!+(r-2)!\bigr)B_{ij}
  =\frac{r!}{r-1}B_{ij}.
\]
The claim follows since this holds for every entry and every $B\in V$.

{\em End of proof of Claim \ref{cl3}}

We now consider the eigenvalues of $L^TL$.
\begin{claim}\label{cl4}
\[
  L^{\mathsf T}L  \text{ has eigenvalues }\frac{K}{r-1}\text{ with multiplicity }d, \quad 0\text{ with multiplicity }q-d.
\]
\end{claim}
The map $L^{\mathsf T}:V\to H$ is injective: let $\lambda=K/(r-1)$.
\[
  L^{\mathsf T}B=0
  \quad\Longrightarrow\quad
  \lambda B=LL^{\mathsf T}B=0.
\]
Thus its image has dimension $d$. For every vector $\bz=L^{\mathsf T}B$ in this image, it follows from Claim \ref{cl3} that 
\[
  L^{\mathsf T}Lz
  =L^{\mathsf T}LL^{\mathsf T}B
  =\lambda L^{\mathsf T}B
  =\lambda \bz.
\]
Thus every nonzero vector in this \(d\)-dimensional subspace is an eigenvector with eigenvalue \(\lambda\).

Also, $L$ is onto $V$, since $B=L\left(\lambda^{-1}L^{\mathsf T}B\right)$ and then $\dim\ker L=q-d$. For $\bz\in\ker L$, we have $L^{\mathsf T}L\bz=0$.

These two subspaces have trivial intersection: if $\bz=L^{\mathsf T}B$ and $L\bz=0$, then $\lambda B=0$, hence $\bz=0$. Their dimensions sum to $q$, so they exhaust $H$. 

Thus $L^TL$ has eigenvalue $K/(r-1)$ with multiplicity $d$, and eigenvalue zero with multiplicity $q-d$.  So, because the eigenvalues of $Q$ are of the form $K-\alpha r \xi$ where $\xi$ is an eigenvalue of $L^TL$,
\[
  \det_H Q=K^q\beta^d, \qquad\text{ where }\beta=\frac{r-s}{s(r-1)}.
\]
Let
\[
  \Lambda=H\cap\mathbb Z^{S_r},
  \qquad
  \mathcal L=V\cap\mathbb Z^{r\times r}.
\]
The map $L:\Lambda\to\mathcal L$ is onto. To see this, let $E_{ij}$ be the $r\times r$ matrix with a one in position $i,j$ and a zero elsewhere. The matrices
\[
  E_{ij}-E_{ir}-E_{rj}+E_{rr}, \qquad 1\le i,j<r,
\]
form an integer basis of the lattice $\mathcal L$, and each is a difference of two permutation matrices.
\begin{claim}\label{cl5}
The sublattice
\[
  \Lambda_s=\{\bz\in\Lambda:L\bz\in s\mathcal L\}
\]
has index
\[
  [\Lambda:\Lambda_s]=|\mathcal L/s\mathcal L|=s^d.
\]
\end{claim}
{\em Proof of Claim \ref{cl5}}
The index $[\Lambda:\Lambda_s]$ is the number of distinct cosets of $\Lambda_s$ in $\Lambda$. We count these by their images under $L$ modulo $s$.

Define
\[
  \varphi:\Lambda\longrightarrow\mathcal L/s\mathcal L,  \qquad  \varphi(\bz)=L\bz+s\mathcal L.
\]
Since $L:\Lambda\to\mathcal L$ is onto, $\varphi$ is onto. Furthermore,
\[
\varphi(\bz)=\varphi(\bz')\iff L(\bz-\bz')\in s\mathcal L\iff\bz-\bz'\in\Lambda_s\iff\bz+\Lambda_s=\bz'+\Lambda_s.
\]
Thus distinct cosets of $\Lambda_s$ correspond exactly to distinct elements of $\mathcal L/s\mathcal L$. Consequently,
\[
  [\Lambda:\Lambda_s]=|\mathcal L/s\mathcal L|.
\]
Recall that $\mathcal L$ consists of integer $r\times r$ matrices with zero row and column sums. Such a matrix is uniquely determined by its upper-left $(r-1)\times(r-1)$ entries: the remaining entries are forced by the zero-sum conditions. 
More explicitly, the matrices
\[
  B_{ij}=E_{ij}-E_{ir}-E_{rj}+E_{rr},
  \qquad 1\le i,j<r,
\]
form an integer basis, and every $M\in\mathcal L$ has the unique representation
\[
  M=\sum_{i,j<r}M_{ij}B_{ij}.
\]
Modulo $s\mathcal L$, each coefficient $M_{ij}$ can independently take one of the $s$ residues $0,\ldots,s-1$. There are $d=(r-1)^2$ coefficients, so 
\[
  |\mathcal L/s\mathcal L|=s^{(r-1)^2}=s^d,
\]
and the claim follows.\\
{\em End of proof of Claim \ref{cl5}}

Recall that the {\em count vector} $c_\sigma$ counts the vertices $v$ for which $\sigma_v=\sigma$. Let $\mathrm{id}$ denote the identity permutation, $\mathrm{id}(i)=i$. Define the vector $\bc^{(0)}$ by $ \bc^{(0)}_{\mathrm{id}}=n$ and $\bc^{(0)}_\sigma=0$ if $\sigma\ne\mathrm{id}$. Its corresponding $N$ matrix (see \eqref{N=}) is $N^{(0)}  =\sum_\sigma \bc^{(0)}_\sigma P_\sigma  =nP_{\mathrm{id}}  =nI$. Every other admissible count vector has the form $\bc=\bc^{(0)}+\bz,\, \bz\in\Lambda_s$, subject to all coordinates of $\bc$ being nonnegative.

Identify $ \Lambda$ with $\left\{\bz\in\mathbb Z^K:\sum_{i=1}^Kz_i=0\right\}$. Let $\be_1,\ldots,\be_K$ be the standard coordinate vectors and set $\bb_i=\be_i-\be_K,1\le i\le q$. These form an integer basis of $\Lambda$.

The Gram matrix $G$ of the basis is
\[
  G_{ij}=\langle \bb_i,\bb_j\rangle  =\langle \be_i-\be_K,\be_j-\be_K\rangle  =\delta_{ij}+1.
\]
Therefore $G=I_q+\mathbf1\mathbf1^{\mathsf T}$. Since $ G\mathbf1=(q+1)\mathbf1=K\mathbf1$ and $Gx=x$ whenever $\mathbf1^{\mathsf T}x=0$, its eigenvalues are $K$ once and 1 with multiplicity $q-1$. Hence $ \det G=K$.

To relate this determinant to volume, choose an orthonormal basis of $H$ and let $A$ be the $q\times q$ matrix whose columns are the coordinates of $\bb_1,\ldots,\bb_q$ in that basis. Let \(\mathcal P\) be the fundamental parallelepiped spanned by \(\bb_1,\ldots,\bb_q\). Then $\operatorname{vol}_H(\mathcal P)=|\det A|$, $G=A^{\mathsf T}A$.

Consequently, $\operatorname{vol}_H(\mathcal P)^2=(\det A)^2=\det(A^{\mathsf T}A)  =\det G=K$. Thus $\operatorname{vol}_H(\mathcal P)=\sqrt K$. The corresponding volume for $\Lambda_s$ is therefore $s^d\sqrt K$.

Ignoring nonnegativity for the moment, the integer count vectors satisfying the total-count and divisibility constraints, after centering and scaling, form the affine lattice
\[
\mathcal Y_n=\frac{\bc^{(0)}-(n/K)\mathbf1}{\sqrt n} +n^{-1/2}\Lambda_s.
\] 
In the central region \(\|\by\|_1\le n^{1/10}\), every corresponding count is positive for sufficiently large \(n\), since
\[
c_\sigma=\frac nK+\sqrt n\,y_\sigma\ge\frac nK-n^{3/5}>0.
\]
\begin{claim}\label{cl6}
\begin{equation}\label{estimate}
n^{-q/2}  \sum_{y\in\mathcal Y_n}e^{-\langle y,Qy\rangle/2}  =  \frac{1}{s^d\sqrt K}  \int_H e^{-\langle \bx,Q\bx\rangle/2}\,d\bx+O(n^{-1/2}).
\end{equation}
\end{claim}
{\em Proof of Claim \ref{cl6}}\\
Put
\[
  h=n^{-1/2},\qquad  \mathcal Y_n=t_n+h\Lambda_s,\qquad  f(x)=e^{-\langle x,Qx\rangle/2}.
\]
Choose a bounded half-open fundamental parallelepiped $P$ of $\Lambda_s$, with
\[
\Xi_P=\operatorname{vol}_H(P)=s^d\sqrt K,  \qquad  D=\sup_{p\in P}\|p\|.
\]
The cells $C_\by=\by+hP,y\in\mathcal Y_n$, partition $H$, apart from boundaries of measure zero. Every cell has volume $h^q\Xi_P$, and $x\in C_\by\Longrightarrow \|\bx-\by\|\le Dh$. 

Let $\lambda>0$ be the smallest eigenvalue of $Q$. Then
\[
  f(\by)\le e^{-\lambda\|\by\|^2/2}.
\]
For $x\in C_\by$, the inequality $\|\bx\|^2\le2\|\by\|^2+2D^2h^2$ gives
\[
  f(\by)  \le  e^{\lambda D^2h^2/2}e^{-\lambda\|\bx\|^2/4} \le C e^{-\lambda\|\bx\|^2/4},
\]
where $C$ is independent of $n$, since $h\le1$.

Therefore, integrating over each cell, and using $\Xi_P=\Theta(1)$,
\[
h^qf(\by)  \le C\int_{C_\by}e^{-\lambda\|\bx\|^2/4}\,d\bx.
\]
If $\|\by\|>R$, then every $\bx\in C_\by$ satisfies $\|\bx\|>R-Dh$. Summing over these cells yields
\[
h^q\sum_{\substack{y\in\mathcal Y_n\\\|\by\|>R}}f(\by)  \le  C\int_{\|\bx\|>R-Dh}e^{-\lambda\|\bx\|^2/4}\,d\bx.
\]
For $T\ge0$,
\[
\int_{\|\bx\|>T}e^{-\lambda\|\bx\|^2/4}\,d\bx  \le  e^{-\lambda T^2/8}  \int_H e^{-\lambda\|\bx\|^2/8}\,d\bx  \le C_1e^{-\lambda T^2/8}.
\]
Consequently, for $R\ge2Dh$,
\begin{equation}\label{hqsum}
h^q\sum_{\substack{y\in\mathcal Y_n\\\|\by\|>R}}f(\by)  \le C_2e^{-\lambda R^2/32}.
\end{equation}

Next let $M$ be the largest eigenvalue of $Q$. Since
\[
  \nabla f(\bz)=-Q\bz\,f(\bz),
\]
we have
\[
\|\nabla f(\bz)\|  \le M\|\bz\|e^{-\lambda\|\bz\|^2/2}.
\]
If $\bx\in C_\by$ and $\bz$ lies on the segment joining $\bx$ to $y$, then $\|\bz-\bx\|\le Dh\le D$. Hence
\[
  \|\nabla f(\bz)\|  \le C_3(1+\|\bx\|)e^{-\lambda\|\bx\|^2/4}.
\]
By the mean value theorem,
\[
  |f(\by)-f(\bx)|  \le  C_4 h(1+\|\bx\|)e^{-\lambda\|\bx\|^2/4}.
\]
The function on the right, after removing the factor $h$, is integrable on $H$. Integrating over all cells therefore gives
\[
 \left|\Xi_Ph^q\sum_{y\in\mathcal Y_n}f(\by)-\int_H f(\bx)\,d\bx\right|\le\sum_{y\in\mathcal Y_n}\int_{C_\by}|f(\by)-f(\bx)|\,d\bx\le  C_4h\int_H(1+\|\bx\|)e^{-\lambda\|\bx\|^2/4}\,d\bx\le C_5h.
\]
Thus, 
\begin{equation}\label{ff1}
 n^{-q/2}\sum_{y\in\mathcal Y_n}f(\by)  =\frac{1}{\Xi_P}\int_H f(\bx)\,d\bx+O(n^{-1/2}).
\end{equation}
Finally, because $\by$ has $K$ coordinates, $ \|\by\|_1\le\sqrt K\,\|\by\|$. Therefore $\|\by\|_1>n^{1/10}$ implies $\|\by\|>n^{1/10}/\sqrt K$. Applying \eqref{hqsum} with this value of $R$ gives
\begin{equation}\label{ff2}
  n^{-q/2}  \sum_{\substack{\by\in\mathcal Y_n\\\|\by\|_1>n^{1/10}}}f(\by)=O\!\left(e^{-\lambda n^{1/5}/(32K)}\right).
\end{equation}
Together with $h=n^{-1/2}$ and $\Xi_P=s^d\sqrt K$, this proves \eqref{estimate}.

{\em End of proof of Claim \ref{cl6}}

Going back to \eqref{muwq} we observe that if $ |\varepsilon_n(\by)|\le C_6n^{-1/5}$ then
\begin{equation}\label{error1}
\left|h^q\sum_{\substack{y\in\mathcal Y_n\\\|\by\|_1\le n^{1/10}}}f(\by)\varepsilon_n(\by)\right|\le  C_6n^{-1/5}h^q\sum_{y\in\mathcal Y_n}f(\by)=O(n^{-1/5}).
\end{equation}
 It follows from \eqref{muwq}, \eqref{estimate}, \eqref{ff1} and \eqref{ff2} that
\begin{equation}\label{final}
\sum_{\substack{\bc\text{ admissible}\\ \|\by(\bc)\|_1\leq n^{1/10}}}\mu_n(\bc)W_n(\bc)=\frac{s^{d/2}K^{K/2}}{s^dK^{1/2}(2\pi)^{q/2}}\int_H e^{-\langle \bx,Q\bx\rangle/2}\,d\bx+O(n^{-1/5}).
\end{equation}
\begin{claim}\label{cl7}
\begin{equation}\label{intQ}
\int_H e^{-\langle \bx,Q\bx\rangle/2}\,d\bx=\frac{(2\pi)^{q/2}}{K^{q/2}\beta^{d/2}},\qquad\text{ where }\beta=\frac{r-s}{s(r-1)}.
\end{equation}
\end{claim}
{\em Proof of Claim \ref{cl7}}

Since \(Q:H\to H\) is symmetric and positive definite, there is an orthonormal basis \(\bv_1,\ldots,\bv_q\) of \(H\) consisting of
eigenvectors: $Q\bv_i=\lambda_i \bv_i,\,\lambda_i>0$. Write $\by=\sum_{i=1}^q z_i \bv_i$. Orthonormality gives
\[
\langle y,Qy\rangle=\sum_{i,j=1}^q z_i z_j\lambda_j\langle \bv_i,\bv_j\rangle=\sum_{i=1}^q\lambda_i z_i^2.
\]
Thus the integral separates into \(q\) one-dimensional integrals:
\[
\int_H e^{-\langle \by,Q\by\rangle/2}\,d\by=\int_{\mathbb R^q}\prod_{i=1}^q e^{-\lambda_i z_i^2/2}\,dz_1\cdots dz_q=\prod_{i=1}^q\int_{\mathbb R}e^{-\lambda_i z_i^2/2}\,dz_i=\prod_{i=1}^q\sqrt{\frac{2\pi}{\lambda_i}}.
\]
Consequently,
\[
\int_H e^{-\langle \by,Q\by\rangle/2}\,d\by=\frac{(2\pi)^{q/2}}{\sqrt{\lambda_1\cdots\lambda_q}}=\frac{(2\pi)^{q/2}}{\sqrt{\det_H Q}}=\frac{(2\pi)^{q/2}}{K^{q/2}\beta^{d/2}}.
\]
{\em End of proof of Claim \ref{cl7}}

Combining \eqref{final} and \eqref{intQ} and using $q=K-1$, we obtain
\begin{equation}\label{ratio}
  \frac{\E_{P_n}Z^2}{(\E_{P_n}Z)^2}  \longrightarrow \frac{s^{d/2}K^{K/2}}{s^dK^{1/2}(2\pi)^{q/2}}\frac{(2\pi)^{q/2}}{K^{q/2}\beta^{d/2}} =s^{-d/2}\beta^{-d/2}=\left(\frac{r-1}{r-s}\right)^{d/2}.
\end{equation}
\end{proof}
\section{Short cycles and mixed factorial moments}\label{shortcycles}
We now proceed with the small subgraph conditioning argument. View a configuration as a bipartite incidence multigraph: there are $n$ left vertices of degree $r$ and $rn/s$ right vertices of degree $s$. A right vertex represents a hyperedge. Let $X_k$ count its cycles of length $2k$. For $k=1$ this means a pair of parallel incidence edges, so $X_1=0$ if and only if the configuration has no loops. 

For fixed $k\ge1$, the usual incidence counting gives
\begin{equation}\label{eq:lambda}
 \lambda_k=\lim_{n\to\infty}\E_{P_n}X_k =\frac{((r-1)(s-1))^k}{2k}.
\end{equation}
For example, temporarily labelling the right vertices, the cycle count has expectation
\[
 \frac{(n)_k(rn/s)_k(r(r-1)s(s-1))^k}{2k(rn)_{2k}} \sim\lambda_k.
\]
Here $(a)_b=a(a-1)\cdots(a-b+1)$. Disjoint cycles give the leading joint factorial moments. Intersecting distinct cycles have at least two independent cycles in one connected incidence component and contribute $O(n^{-1})$. It follows that, for each fixed $m$,
\begin{equation}\label{eq:poisson}
 (X_1,\ldots,X_m)\ \xrightarrow{\ d\ }\ (\xi_1,\ldots,\xi_m), \qquad \xi_k\ \text{independent }\operatorname{Poisson}(\lambda_k).
\end{equation}
Under $Q_n$, right vertices have matching colours, and consecutive right vertices on a cycle must have distinct colours. The number of such colourings of a cycle of length $k$ is
\[
 C_k(r)=(r-1)^k+(-1)^k(r-1).
\]
This formula also gives $C_1(r)=0$ and $C_2(r)=r(r-1)$. For each admissible colouring, incidence counting contributes $(s-1)^k/(2k)$ asymptotically. Thus
\begin{equation}\label{eq:mu}
 \mu_k:=\lim_{n\to\infty}\E_{Q_n}X_k =\frac{(s-1)^k}{2k} \big((r-1)^k+(-1)^k(r-1)\big).
\end{equation}
The same disjoint-cycle argument gives the joint factorial moments under $Q_n$. Equivalently, for every fixed list of nonnegative integers $\bb_1,\ldots,\bb_m$,
\begin{equation}\label{eq:mixed}
 \frac{\E_{P_n}(Z\prod_{k=1}^m(X_k)_{\bb_k})}{\E_{P_n}Z} \longrightarrow \prod_{k=1}^m\mu_k^{\bb_k}.
\end{equation}
Define
\begin{equation}\label{eq:delta}
 \delta_k=\frac{\mu_k}{\lambda_k}-1 =\frac{(-1)^k}{(r-1)^{k-1}}.
\end{equation}
In particular, $\delta_1=-1$, as expected since $Q_n$ has no loops. Since $r>s$,
\begin{equation}\label{eq:cyclesum}
\sum_{k\ge1}\lambda_k\delta_k^2=\frac{(r-1)^2}{2}\sum_{k\ge1}\frac1k\left(\frac{s-1}{r-1}\right)^k =-\frac{(r-1)^2}{2} \log\left(1-\frac{s-1}{r-1}\right) <\infty.
 \end{equation}
Comparing \eqref{ratio} and \eqref{eq:cyclesum} we see that we can write
\begin{equation}\label{limit}
 R_n=\frac{\E_{P_n}Z^2}{(\E_{P_n}Z)^2}  \longrightarrow \exp\left\{\sum_{k\ge1}\lambda_k\delta_k^2\right\}=e^{\lambda_1+S},
\end{equation}
since $\delta_1=-1$ and if we let $S$ denote the sum from 2 to $\infty$.
\section{Small-subgraph conditioning and mutual contiguity}\label{sec:contiguity}
Throughout this section, $s\ge3$. We work with the configuration model conditioned to have no loops:
\[
 P_n^\circ=P_n(\,\cdot\mid X_1=0),\qquad p_n=P_n(X_1=0)\longrightarrow e^{-\lambda_1}>0, \qquad \lambda_1=\frac{(r-1)(s-1)}2.
\]
The conditioning is necessary: even when $s\ge3$, the unconditioned configuration model has $X_1\xrightarrow{\mathrm d}\operatorname{Poisson}(\lambda_1)$, so $P_n(X_1=0)$ does not tend to one. The assumption $s\ge3$ will instead make repeated whole hyperedges negligible at the end of the proof.

Now let 
\[
L_n(F)=\frac{Q_n(F)}{P_n^\circ(F)}=p_n\frac{Z(F)}{\E_{P_n}Z},
\]
and \(L_n=0\) on configurations with loops. 

Then for any event $\mathcal A$, we have
\begin{equation}\label{event}
\mathbb E_{P_n^\circ}[L_n\mathbf1_{\mathcal A}]=\sum_F L_n(F)\mathbf1_{\mathcal A}(F)P_n^\circ(F)=\sum_{F\in\mathcal A}L_n(F)P_n^\circ(F)=\sum_{F\in\mathcal A}Q_n(F)=Q_n(\mathcal A).
\end{equation}
Furthermore,
\[
\mathbb E_{P_n^\circ}L_n^2=\frac{\mathbb E_{P_n^\circ}Z^2}{(\mathbb E_{P_n^\circ}Z)^2}=\frac{\mathbb E_{P_n}Z^2/p_n}{(\mathbb E_{P_n}Z/p_n)^2}=p_n\frac{\mathbb E_{P_n}Z^2}{(\mathbb E_{P_n}Z)^2} \longrightarrow e^{-\lambda_1} \left(\frac{r-1}{r-s}\right)^{(r-1)^2/2}.
\]
For every fixed $m\ge2$, under $P_n^\circ$,
\[
 (X_2,\ldots,X_m)\xrightarrow{\mathrm d} (\xi_2,\ldots,\xi_m),\qquad
 \xi_k\text{ independent }\operatorname{Poisson}(\lambda_k).
\]
Indeed, joint convergence in \eqref{eq:poisson}, together with $p_n\to e^{-\lambda_1}$ and independence of the limiting variables, gives this conditional convergence. Under $Q_n$, the factorial moments from the preceding section give independent Poisson limits with means  $\mu_k=\lambda_k(1+\delta_k)$. Equivalently, for fixed nonnegative integers $\bb_2,\ldots,\bb_m$,
\begin{claim}\label{cl8}
\begin{equation}\label{tosee}
 \E_{P_n^\circ}\left[L_n\prod_{k=2}^m(X_k)_{\bb_k}\right] =\frac{\E_{P_n}\left[Z\prod_{k=2}^m(X_k)_{\bb_k}\right]}{\E_{P_n}Z} \longrightarrow\prod_{k=2}^m\mu_k^{\bb_k}.
\end{equation}
\end{claim}
{\em Proof of Claim \ref{cl8}}\\
Put $Y_n=\prod_{k=2}^m(X_k)_{\bb_k}$. Since \(Z=0\) whenever \(X_1>0\), the product \(ZY_n\) is also zero there. Consequently,
\[
\mathbb E_{P_n^\circ}[ZY_n]=\frac{\mathbb E_{P_n}[ZY_n]}{p_n},\qquad\mathbb E_{P_n^\circ}Z=\frac{\mathbb E_{P_n}Z}{p_n}.
\]
Using \(L_n=Z/\mathbb E_{P_n^\circ}Z\), we therefore obtain
\[
\mathbb E_{P_n^\circ}[L_nY_n]=\frac{\mathbb E_{P_n^\circ}[ZY_n]}{\mathbb E_{P_n^\circ}Z}=\frac{\mathbb E_{P_n}[ZY_n]/p_n}{\mathbb E_{P_n}Z/p_n}=\frac{\mathbb E_{P_n}[ZY_n]}{\mathbb E_{P_n}Z}.
\]
For the limit, the definition of \(Q_n\) gives $\frac{\mathbb E_{P_n}[ZY_n]}{\mathbb E_{P_n}Z}=\mathbb E_{Q_n}Y_n$. The cycle-counting calculation in Section~\ref{shortcycles} proves
\[
\mathbb E_{Q_n}\prod_{k=2}^m(X_k)_{\bb_k}\longrightarrow\prod_{k=2}^m\mu_k^{\bb_k}.
\]
Indeed, \((X_k)_{\bb_k}\) counts ordered choices of \(\bb_k\) distinct cycles of length \(2k\). Collections in which all the chosen cycles are disjoint contribute \(\prod_{k=2}^m\mu_k^{\bb_k}+o(1)\), while collections containing intersecting cycles contribute \(o(1)\). This completes the proof of \eqref{tosee}.\\
{\em End of proof of Claim \ref{cl8}}

Define
\[
 W_m=\prod_{k=2}^m(1+\delta_k)^{\xi_k}e^{-\lambda_k\delta_k}.
\]
Now let
\[
S_m:=\sum_{k=2}^m\lambda_k\delta_k^2,\quad  S:=\sum_{k=2}^{\infty}\lambda_k\delta_k^2.
\]
Then
\[
 \E W_m=1,\qquad \E W_m^2=e^{S_m}.
\]
Thus $(W_m)$ is a martingale, bounded in $L^2$, with respect to the sigma-fields generated by $\xi_2,\ldots,\xi_m$. It converges almost surely and in $L^2$ to a random variable $W$, with
\begin{equation}\label{eq:Wmoments}
 \E W=1,\qquad \E W^2=e^S.
\end{equation}
Let $\boldsymbol X_{n,m}=(X_2,\ldots,X_m)$ and put
\[
 U_{n,m}=\E_{P_n^\circ}[L_n\mid\boldsymbol X_{n,m}].
\]
Fix \(m\ge2\) and a vector \(\boldsymbol x=(x_2,\ldots,x_m)\) of nonnegative integers. Write
\[
\mathcal A_{n,\boldsymbol x}=\{\boldsymbol X_{n,m}=\boldsymbol x\}=\{X_2=x_2,\ldots,X_m=x_m\}.
\]
Then
\[
\mathbb E_{P_n^\circ}[L_n\mid\mathcal A_{n,\boldsymbol x}]=\frac{\mathbb E_{P_n^\circ}[L_n\mathbf1_{\mathcal A_{n,\boldsymbol x}}]}{P_n^\circ(\mathcal A_{n,\boldsymbol x})}=\frac{Q_n(\mathcal A_{n,\boldsymbol x})}{P_n^\circ(\mathcal A_{n,\boldsymbol x})},
\]
after using \eqref{event}.

Under \(P_n^\circ\), the vector \((X_2,\ldots,X_m)\) converges to independent Poisson variables with means\\ \(\lambda_2,\ldots,\lambda_m\). Therefore
\[
P_n^\circ(\mathcal A_{n,\boldsymbol x})\longrightarrow\prod_{k=2}^m e^{-\lambda_k}      \frac{\lambda_k^{x_k}}{x_k!}.
\]
Under \(Q_n\), the corresponding limiting means are \(\mu_2,\ldots,\mu_m\), so
\[
Q_n(\mathcal A_{n,\boldsymbol x})\longrightarrow\prod_{k=2}^m e^{-\mu_k} \frac{\mu_k^{x_k}}{x_k!}.
\]
The limiting denominator is strictly positive. Thus 
\begin{equation}\label{glen1}
 \E_{P_n^\circ}[L_n\mid\boldsymbol X_{n,m}=\boldsymbol x] =\frac{Q_n(\mathcal A_{n,\boldsymbol x})}{P_n^\circ(\mathcal A_{n,\boldsymbol x})} \longrightarrow \prod_{k=2}^m(1+\delta_k)^{x_k}e^{-\lambda_k\delta_k}.
\end{equation}

Keep \(m\) fixed. Write
\[
g_n(\boldsymbol x)=\mathbb E_{P_n^\circ}  [L_n\mid\boldsymbol X_{n,m}=\boldsymbol x],\qquad g(\boldsymbol x)=\prod_{k=2}^m (1+\delta_k)^{x_k}e^{-\lambda_k\delta_k}.
\]
Then
\begin{equation}\label{glen2}
U_{n,m}=g_n(\boldsymbol X_{n,m}),\qquad W_m=g(\boldsymbol\xi_m),
\end{equation}
where \(\boldsymbol\xi_m=(\xi_2,\ldots,\xi_m)\). We have proved that
\[
g_n(\boldsymbol x)\longrightarrow g(\boldsymbol x)\quad\text{for each fixed }\boldsymbol x,\qquad\boldsymbol X_{n,m}\xrightarrow{\mathrm d}\boldsymbol\xi_m.
\]
Choose \(R\) so that the finite box $B_R=\{0,1,\ldots,R\}^{m-1}$ satisfies
\[
\mathbb P(\boldsymbol\xi_m\notin B_R)<\varepsilon,\qquad P_n^\circ(\boldsymbol X_{n,m}\notin B_R)<\varepsilon
\]
for all sufficiently large \(n\). Indeed, first choose the limiting probability to be less than \(\varepsilon/2\), and then use convergence of the probabilities of the finite set \(B_R\).

Let \(\varphi:\mathbb R\to\mathbb R\) be bounded and continuous, and put
\(M=\sup_t|\varphi(t)|\). For each \(\boldsymbol x\in B_R\),
\[
\varphi(g_n(\boldsymbol x))P_n^\circ(\boldsymbol X_{n,m}=\boldsymbol x)\longrightarrow\varphi(g(\boldsymbol x))\mathbb P(\boldsymbol\xi_m=\boldsymbol x).
\]
Since \(B_R\) is finite, we may sum these convergences:
\[
\mathbb E_{P_n^\circ} \left[\varphi(U_{n,m})\mathbf1_{\{\boldsymbol X_{n,m}\in B_R\}}\right]\longrightarrow\mathbb E\left[\varphi(W_m)\mathbf1_{\{\boldsymbol\xi_m\in B_R\}}\right].
\]
Each omitted tail has absolute value at most \(M\varepsilon\). Therefore
\[
\limsup_{n\to\infty}\left|\mathbb E_{P_n^\circ}\varphi(U_{n,m})-\mathbb E\varphi(W_m)\right|\le 2M\varepsilon.
\]
Letting \(\varepsilon\downarrow0\) proves convergence of expectations for every bounded continuous \(\varphi\), which is equivalent to
\[
  U_{n,m}\xrightarrow{\mathrm d}W_m.
\]
Summing over any finite set of vectors and then increasing that set also gives
\[
 \liminf_{n\to\infty}\E_{P_n^\circ}U_{n,m}^2 \ge\E W_m^2=e^{S_m}.
\]

The residual \(L_n-U_{n,m}\) has conditional mean zero: because \(U_{n,m}\) is measurable with respect to \(\boldsymbol X_{n,m}\), it can be taken outside this conditional expectation. Then, from \eqref{glen1}, \eqref{glen2},
\[
\mathbb E_{P_n^\circ} [L_n-U_{n,m}\mid\boldsymbol X_{n,m}]=U_{n,m}-U_{n,m}=0.
\]
Define
\[
 \mathcal F_{n,m}:=\sigma(X_2,\ldots,X_m), \qquad
 U_{n,m}:=\E_{P_n^\circ}[L_n\mid\mathcal F_{n,m}].
\]
 Using the tower property, we obtain
\[
\mathbb E_{P_n^\circ}[(L_n-U_{n,m})U_{n,m}]=\mathbb E_{P_n^\circ}\left[
 U_{n,m}\, \mathbb E_{P_n^\circ} [L_n-U_{n,m}\mid\mathcal F_{n,m}]\right]=0.
\]
Consequently, 
\[
 \E_{P_n^\circ}(L_n-U_{n,m})^2 =\E_{P_n^\circ}L_n^2-\E_{P_n^\circ}U_{n,m}^2.
\]
Therefore
\[
 \limsup_{n\to\infty}
 \E_{P_n^\circ}(L_n-U_{n,m})^2\le e^S-e^{S_m}
 \longrightarrow0\qquad(m\to\infty).
\]
Together with $U_{n,m}\xrightarrow{\mathrm d}W_m$ for fixed $m$ and $W_m\to W$ in $L^2$, this proves

\begin{claim}\label{cl9}
\begin{equation}\label{eq:W}
 L_n\xrightarrow{\mathrm d}W =\prod_{k=2}^{\infty}(1+\delta_k)^{\xi_k}e^{-\lambda_k\delta_k}.
\end{equation}
\end{claim}
{\em Proof of Claim \ref{cl9}}\\
Let \(\varphi:\mathbb R\to\mathbb R\) be bounded and Lipschitz,
with Lipschitz constant \(C_\varphi\). 
The triangle inequality gives
\begin{multline*}
\left|\mathbb E_{P_n^\circ}\varphi(L_n)-\mathbb E\varphi(W)\right|\\
\le\left|\mathbb E_{P_n^\circ}\varphi(L_n)-\mathbb E_{P_n^\circ}\varphi(U_{n,m})\right|+\left|\mathbb E_{P_n^\circ}\varphi(U_{n,m})-\mathbb E\varphi(W_m)\right|+\left|\mathbb E\varphi(W_m)-\mathbb E\varphi(W)\right|.
\end{multline*}

For the first term, the Lipschitz property followed by Cauchy--Schwarz gives
\[
\left|\mathbb E_{P_n^\circ}\varphi(L_n)-\mathbb E_{P_n^\circ}\varphi(U_{n,m})\right|\le C_\varphi\mathbb E_{P_n^\circ}|L_n-U_{n,m}|\le C_\varphi\left(\mathbb E_{P_n^\circ}(L_n-U_{n,m})^2\right)^{1/2}.
\]
The same argument bounds the third term by $C_\varphi\left(\mathbb E(W_m-W)^2\right)^{1/2}$. Now fix \(m\) and let \(n\to\infty\). The middle term tends to zero because \(U_{n,m}\xrightarrow{\mathrm d}W_m\) and
\(\varphi\) is bounded and continuous. Therefore
\[
\limsup_{n\to\infty}\left|\mathbb E_{P_n^\circ}\varphi(L_n)-\mathbb E\varphi(W)\right|\le C_\varphi\sqrt{e^S-e^{S_m}}+C_\varphi\left(\mathbb E(W_m-W)^2\right)^{1/2}.
\]
Finally let \(m\to\infty\). The first term tends to zero because \(S_m\to S\), and the second tends to zero because \(W_m\to W\) in \(L^2\). Hence $\mathbb E_{P_n^\circ}\varphi(L_n)\longrightarrow\mathbb E\varphi(W)$
for every bounded Lipschitz function \(\varphi\). These test functions characterize convergence in distribution. So $L_n\xrightarrow{\mathrm d}W$.\\
{\em End of proof of Claim \ref{cl9}}\\
The limit is strictly positive almost surely. Indeed, $\delta_k>-1$ for every $k\ge2$, and
\[
\log W_m=\sum_{k=2}^m(\xi_k-\lambda_k)\log(1+\delta_k)+\sum_{k=2}^m\lambda_k\bigl(\log(1+\delta_k)-\delta_k\bigr).
\]
Since $\delta_k\to0$ and $S<\infty$,
\[
 \sum_{k\ge2}\lambda_k\log^2(1+\delta_k)<\infty,\qquad \sum_{k\ge2}\lambda_k|\log(1+\delta_k)-\delta_k|<\infty.
\]
The first series in the expression for $\log W_m$ converges almost surely, because its independent centred summands have summable variances. The second converges absolutely. Thus $\log W_m$ has a finite real limit, and $W>0$ almost surely.
\paragraph{Mutual contiguity.}
Let $(\mathcal A_n)$ be any sequence of configuration events. If $P_n^\circ(\mathcal A_n)\to0$, then
\[
 Q_n(\mathcal A_n)=\E_{P_n^\circ}[L_n\mathbf1_{\mathcal A_n}]\le\bigl(\E_{P_n^\circ}L_n^2\bigr)^{1/2}P_n^\circ(\mathcal A_n)^{1/2}\longrightarrow0.
\]
Conversely, for every $\varepsilon>0$,
\[
 P_n^\circ(\mathcal A_n)\le P_n^\circ(L_n\le\varepsilon)+\varepsilon^{-1}Q_n(\mathcal A_n).
\]
If $Q_n(\mathcal A_n)\to0$, then equation \eqref{eq:W} implies that $\limsup_{n\to\infty}P_n^\circ(\mathcal A_n)\le\Pp(W\le\varepsilon)$. Letting $\varepsilon\downarrow0$ and using $W>0$ almost surely shows that $P_n^\circ(\mathcal A_n)\to0$. Thus $P_n^\circ$ and $Q_n$ are mutually contiguous.

\paragraph{Projection to simple hypergraphs.}
Let $\mathcal S_n$ be the event that the projected configuration is a simple hypergraph. Under $Q_n$, write its independent perfect matchings as $M_1,\ldots,M_r$. There are no loops, and any repeated hyperedge must occur in two different matchings. For a fixed $s$-set $e$, symmetry gives
\[
 \Pp(e\in M_a)=\frac{n/s}{\binom ns}.
\]
Independence of the matchings therefore implies
\[
  \E_{Q_n}\sum_{1\le a<b\le r}|M_a\cap M_b|=\binom r2\frac{(n/s)^2}{\binom ns}=O(n^{2-s})=o(1),
\]
where we use $s\ge3$. Hence $Q_n(\mathcal S_n)\to1$. Mutual contiguity gives $P_n^\circ(\mathcal S_n)\to1$ as well.

Under $P_n^\circ(\,\cdot\mid\mathcal S_n)$, projection gives the uniform simple $r$-regular $s$-uniform hypergraph: each such labelled hypergraph has exactly $(r!)^n$ configuration preimages, obtained by assigning its incident edges to the distinguishable points at each vertex. Conditioning on an event of probability $1-o(1)$ preserves contiguity, as does applying the projection map. The projection of $Q_n$ is $U_{n,r,s}$, so this proves Theorem~\ref{thm:main}.

\appendix
\section{Proof of \eqref{eq:BC}}\label{sec:BC}
For a function $f:S_r\to \reals_+=[0,\infty)$ let 
\begin{align*}
\E_u(f)&=\frac{1}{r!}\sum_{\sigma\in S_r}f(\sigma).\\
\Ent_u(f)&=\frac1{r!}\sum_{\sigma\in S_r}f(\sigma)\log f(\sigma)-\E_uf\log \E_uf.
\end{align*}
For $i\in[r]$, define the coordinate function
\[
  X_i:S_r\longrightarrow[r],\qquad X_i(\rho)=\rho(i).
\]
For a probability distribution $\pi$ on $S_r$ and $f:S_r\to\reals_+$,  define
\[
  h_i^\pi(f,\tau)  :=\bigl(\mathbb E_\pi[f\mid X_i]\bigr)(\tau)  =  \frac{   \sum_{{\rho\in S_r,\rho(i)=\tau(i)}} \pi(\rho)f(\rho)  }{\sum_{{\rho\in S_r,\rho(i)=\tau(i)}}\pi(\rho)},
\]
whenever the denominator is positive. In particular,
\[
 h_i^u(f,\tau)= \frac{1}{(r-1)!}\sum_{\substack{\rho\in S_r\\\rho(i)=\tau(i)}}f(\rho).
\]
Bristiel and Caputo, equation~(4.13), prove that for every $f:S_r\to\reals_+$,
\begin{equation}\label{BCf}
\sum_{i=1}^r\operatorname{Ent}_u\!\left(\mathbb E_u[f\mid X_i]\right)\leq
  \frac{r\log r}{\log(r!)}\operatorname{Ent}_u(f).
\end{equation}
This uses their notation with $\mu=u$, $n=r$, $x=i$, and $\sigma_x=X_i$. For a proof of \eqref{BCf}, see the next section.

We apply \eqref{BCf} with $f\gets g=r!p$ to obtain 
\begin{equation}\label{BCg}
\sum_{i=1}^r  \operatorname{Ent}_u(\E_u(g\mid X_i))\leq  \frac{r\log r}{\log(r!)}\operatorname{Ent}_u(g).
\end{equation}
What remains is to show that
\begin{equation}\label{toshow}
\sum_{i=1}^r  \operatorname{Ent}_u \bigl(\mathbb E_u[g\mid X_i]\bigr)=B(\bp)  \quad\text{and}\quad\operatorname{Ent}_u(g)=D(\bp\Vert \bu).
\end{equation}
Now,
\[
\E_u g=\frac{1}{r!}\sum_{\sigma\in S_r}g(\sigma)=\sum_{\sigma\in S_r}p_\sigma=1,
\]
and so 
\begin{equation}
 \Ent_u(g) =\frac1{r!}\sum_{\sigma\in S_r}r!p_\sigma\log(r!p_\sigma) =D(\bp\Vert \bu).
 \label{eq:full}
\end{equation}
This deals with the second term in \eqref{toshow}.

For fixed $i,j$, there are $(r-1)!$ permutations satisfying $X_i=j$. Under $u$, the conditional distribution on these permutations is uniform. Thus
\[
h_i(\tau):=h_i^u(g,\tau)
  =\bigl(\mathbb E_u[g\mid X_i]\bigr)(\tau)
  =r a_{i,\tau(i)},\qquad \tau\in S_r.
\]
Consequently, since $X_i$ is uniform on $[r]$ under $u$ and $\E_u h_i=1$, 
\[
 \Ent_u(h_i) =\frac1r\sum_{j=1}^r ra_{ij}\log(ra_{ij}) =\sum_{j=1}^r a_{ij}\log(ra_{ij}).
\]
Summing over $i$ gives
\begin{equation}
 \sum_{i=1}^r\Ent_u(h_i)=B(\bp).\label{eq:marginals}
\end{equation}
This deals with the first term in \eqref{toshow}.
\section{Proof of \eqref{BCf}}\label{BC1}
Write \(H\) for Shannon entropy:
\[
H(\bp)=-\sum_jp_j\log p_j.
\]
We first prove that every random permutation \(X=(X_1,\ldots,X_r)\), with an arbitrary distribution on \(S_r\), satisfies
\begin{equation}\label{1}
H(X)\le \frac{\log(r!)}{r\log r}\sum_{i=1}^r H(X_i).
\end{equation}

We need the following estimate.
\begin{claim}\label{cl10}
 If \(\bp=(p_1,\ldots,p_r)\) is a probability vector and
\[
q_j=\frac{1-p_j}{r-1},
\qquad
\theta_r=\frac{\log(r/(r-1))}{\log r},
\]
then
\begin{equation}\label{2}
H(q)\ge \log(r-1)+\theta_r H(\bp).
\end{equation}
\end{claim}
{\em Proof of Claim \ref{cl10}}\\
For \(r=2\), it is equality. For \(r\ge3\), put \(\gamma=(r-1)\theta_r\) and
\[
\phi(t)=-(1-t)\log(1-t)+\gamma t\log t.
\]
The desired inequality is equivalent to
\begin{equation}\label{3}
\sum_{j=1}^r\phi(p_j)\ge0,
\end{equation}
since its left-hand side divided by \(r-1\) equals \(H(q)-\log(r-1)-\theta_rH(\bp)\). Notice that
\[
\phi(0)=\phi(1)=\phi(1/r)=0.
\]
Minimize the left-hand side of \eqref{3} over the probability simplex. A boundary point other than a vertex cannot be a minimizer: moving mass \(\varepsilon\) from a coordinate in \((0,1)\) to a zero coordinate changes the sum by
\[
\gamma\varepsilon\log\varepsilon+O(\varepsilon)<0.
\]
At an interior minimizer, all the values \(\phi'(\bp_j)\) are equal. Moreover,
\[
\phi''(t)=\frac{\gamma}{t}-\frac1{1-t},
\]
so \(\phi'\) first increases and then decreases. Consequently, the coordinates have at most two distinct values. If there are two values, the larger one lies where \(\phi''<0\) and cannot occur twice: perturbing two such coordinates by
\(+\varepsilon\) and \(-\varepsilon\) would decrease the sum. Thus a nonuniform interior minimizer must have the form
\[
\left(t,\frac{1-t}{r-1},\ldots,\frac{1-t}{r-1}\right), \qquad \frac1r<t<1.
\]

Define
\[
F(t)=\phi(t)+(r-1)\phi\!\left(\frac{1-t}{r-1}\right).
\]
We have
\[
F(1/r)=F(1)=0,\qquad F'(1/r)=0,
\]
and direct differentiation gives
\[
F''(t)=
\frac{\gamma(r-2)-(r-1-\gamma)t}{t(1-t)(r-2+t)}.
\]
Since \(\gamma<r-1\), the numerator decreases with \(t\). Hence \(F''\) changes sign at most once, from positive to negative. Starting from \(F'(1/r)=0\), this means that \(F\) can first increase and then decrease, so its minimum is at an endpoint. Therefore \(F(t)\ge0\). The uniform vector also gives zero, which proves \eqref{3}, and hence \eqref{2}.\\
{\em End of proof of Claim \ref{cl10}}

We now prove \eqref{1} by induction on \(r\). The case \(r=2\) is immediate, since either coordinate determines the permutation.

For \(r\ge3\), fix \(i\), choose \(J\) uniformly from \([r]\setminus\{i\}\), independently of \(X\), and put \(Y=X_J\). Conditional on \(X_i=a\), the variable \(Y\) is uniform on \([r]\setminus\{a\}\): each remaining label occurs at exactly one of the \(r-1\) eligible positions. Thus, if \(\bp\) is the distribution of \(X_i\), the distribution of \(Y\) is the vector \(q\) in \eqref{2}, and
\[
H(Y\mid X_i)=\log(r-1).
\]
The entropy chain rule and \eqref{2} give
\[
H(X_i\mid Y)=H(X_i)+H(Y\mid X_i)-H(Y)\le (1-\theta_r)H(X_i)=\frac{\log(r-1)}{\log r}H(X_i).
\]
Since conditioning reduces entropy and \(J\) is independent of \(X\),
\[
\frac1{r-1}\sum_{j\ne i}H(X_i\mid X_j)=H(X_i\mid Y,J)\le H(X_i\mid Y).
\]
Consequently,
\begin{equation}\label{4}
\sum_{j\ne i}H(X_i\mid X_j)\le (r-1)\frac{\log(r-1)}{\log r}H(X_i).
\end{equation}
Conditional on any value of \(X_j\) having positive probability, the remaining coordinates form a random permutation of \(r-1\) labels. Applying the induction hypothesis to each such conditional distribution and averaging gives
\[
H(X\mid X_j)\le \frac{\log((r-1)!)}{(r-1)\log(r-1)}\sum_{i\ne j}H(X_i\mid X_j).
\]
Put \(A=\sum_iH(X_i)\). Summing over \(j\), using the chain rule and then \eqref{4}, yields
\[
rH(X)-A=\sum_{j=1}^rH(X\mid X_j)\le\frac{\log((r-1)!)}{(r-1)\log(r-1)} \sum_{i=1}^r\sum_{j\ne i}H(X_i\mid X_j)\le \frac{\log((r-1)!)}{\log r}\,A.
\]
Rearranging proves \eqref{1}. Finally, we translate \eqref{1} into \eqref{BCf}. Suppose that \(\E_u f=m>0\); the case \(\E_u f=0\) is trivial. Give \(X\) the distribution
\begin{equation}\label{defX}
\mathbb P(X=\sigma)=\frac{f(\sigma)}{r!m},\qquad \sigma\in S_r.
\end{equation}
Then
\begin{claim}\label{cl11}
\[
\Ent_u(f)=m(\log(r!)-H(X)).
\]
\end{claim}
{\em Proof of Claim \ref{cl11}}\\
The Shannon entropy of \(X\) is
\[
\begin{aligned}
H(X)&=-\sum_{\sigma\in S_r}p_\sigma\log p_\sigma=-\frac1{r!m}\sum_{\sigma\in S_r}f(\sigma)\log\!\left(\frac{f(\sigma)}{r!m}\right)\\
&=-\frac1{r!m}\sum_{\sigma\in S_r}f(\sigma)\log f(\sigma)+\frac{\log(r!m)}{r!m}\sum_{\sigma\in S_r}f(\sigma)=\log(r!m)-\frac{\mathbb E_u[f\log f]}m.
\end{aligned}
\]
On the other hand, 
\[
\operatorname{Ent}_u(f)=\mathbb E_u[f\log f] -(\mathbb E_u f)\log(\mathbb E_u f)=\mathbb E_u[f\log f]-m\log m.
\]
Combining these equations gives
\[
\operatorname{Ent}_u(f)=m(\log(r!)-H(X)).
\]
{\em End of proof of Claim \ref{cl11}}\\ 
Write \(X_i=X(i)\) and \(p_i(j)=\mathbb P(X_i=j)\). Let \(\sigma\) denote the identity random permutation on the probability space \((S_r,\bu)\), so that \(\sigma\) is uniformly distributed on \(S_r\). Then, for every fixed \(\tau\in S_r\),
\begin{align*}
 h_i(\tau)&=\bigl(\E_u[f\mid\sigma(i)]\bigr)(\tau)=\frac1{(r-1)!}\sum_{\substack{\rho\in S_r\\\rho(i)=\tau(i)}}f(\rho)=rm\,p_i(\tau(i)).\\
\E_u h_i&=\frac1r\sum_{j=1}^r rm\,p_i(j)=m.
\end{align*}
Consequently,
\[
\begin{aligned}
\operatorname{Ent}_u(h_i)&=\frac1r\sum_{j=1}^rrm\,p_i(j)\log\bigl(rm\,p_i(j)\bigr)-m\log m\\
&=m\log(rm)+m\sum_{j=1}^r p_i(j)\log p_i(j)-m\log m\\
&=m\bigl(\log r-H(X_i)\bigr).
\end{aligned}
\]
Thus \eqref{1} implies
\[
\begin{aligned}
\sum_{i=1}^r\Ent_u\!\left(\E_u[f\mid\sigma(i)]\right)
&=m\left(r\log r-\sum_{i=1}^rH(X_i)\right)\le m\left(r\log r-\frac{r\log r}{\log(r!)} H(X)\right)\\
&=\frac{r\log r}{\log(r!)}\Ent_u(f),
\end{aligned}
\]
as required.

The constant is sharp: for
\(f=r!\mathbf1_{\{\sigma=\sigma_0\}}\), the right-hand entropy
is \(\log(r!)\) and each of the \(r\) marginal entropies is
\(\log r\).

\end{document}